\documentclass[
    ,final            
  ]
  {aipproc}

\layoutstyle{6x9}
\usepackage{amssymb,amsmath}

\begin{document}

\title{Sequences of Knots and Their Limits}


\author{Pedro Lopes}{
  address={Department of Mathematics, Instituto Superior T\'ecnico, Technical University of Lisbon,
Av.~Rovisco Pais, 1049-001 Lisbon Portugal}}

\begin{abstract}
Hyperfinite knots, or limits of equivalence classes of knots induced by a knot invariant taking values in a metric space, were introduced in a previous article by the author. In this article, we present new examples of hyperfinite knots stemming from sequences of torus knots.
\end{abstract}

\maketitle


\section{Introduction}

 Consider a sequence of knots of increasing crossing number, where crossing number of a knot stands for the least number of crossings it takes to draw a diagram of the corresponding knot. Does it make sense to discuss the limit of such a sequence? It does, to a certain extent. Assume there is a knot invariant, say $f$, which takes on values in a complete metric space, call it $\mathbf{M}$. Define two knots to be related if and only if they have the same $f$-invariant. Let ${\cal K}\sb{f}$ be the set of the equivalence classes of the knots under this (equivalence) relation. There is then an $f\sp{\sim}$ map from ${\cal K}\sb{f}$ to $\mathbf{M}$ which is obviously injective. Hence, ${\cal K}\sb{f}$ can be regarded as a metric subspace of $\mathbf{M}$. Let then $\overline{{\cal K}\sb{f}}$ be the closure of ${\cal K}\sb{f}$ in the topology of $\mathbf{M}$. In this way, whenever the sequence of $f$-invariants of a given sequence of knots converges in $\mathbf{M}$, then the corresponding sequence of equivalence classes of the knots converges in $\overline{{\cal K}\sb{f}}$. We call the latter limit a hyperfinite knot (\cite{Lopes}).

\bigbreak

One way of implementing these ideas is via the CJKLS invariant of knots (\cite{Lopes, jsCetal}). There is a CJKLS invariant for each choice of labeling quandle $X$, abelian group $A$, and $2$-co-cycle $\phi$ in $H\sp{2}(X, A)$. Furthermore, the CJKLS invariant is a sum-over-states, where the states are the colorings of the knot under study by the labeling quandle $X$, and the for each coloring the $\phi$'s are evaluated at each crossing of the knot diagram and then multiplied. Since the CJKLS invariant takes on values on the group-algebra of the abelian group $A$, we embed it in the euclidean space $\mathbb{R}\sp{|A|}$, thereby obtaining an invariant taking on values in a complete metric space. We then take logarithms of each coordinate, obtaining what we call the free energy invariant, and then by dividing by the crossing number we obtain what we call the free energy per crossing. We were clearly inspired by Thermodynamics and by Statistical Mechanics of Exactly Solved Models (\cite{Callen, Baxter}).

\bigbreak

\section{Hyperfinite knots whose free energies per crossing equal the null vector}

Specifically, a quandle is an algebraic structure i.e. a set, $X$, endowed with a binary operation, $\ast $, such that for all $a, b, c \in X$:
\begin{itemize}
\item $a\ast a = a$;
\item there is a unique $x \in X$ s.t. $x\ast b = a$;
\item $(a\ast b)\ast c = (a\ast c)\ast (b\ast c)$
\end{itemize}
Any group with conjugation as the $\ast$ operation is a quandle. The ring of Laurent polynomials on the variable $T$, with integer coefficients, $\mathbb{Z}[T, T\sp{-1}]$, endowed with the operation $a\ast b = Ta + (1-T)b$ is also a quandle. Quotients obtained by factoring this quandle by certain ideals in the ring of Laurent polynomials are finite quandles called Alexander quandles. Any knot (oriented) diagram provides a quandle by regarding the arcs of the diagram as generators and reading relations at the crossings of the sort under-arc $\ast $ over-arc equals under-arc where the co-orientation at the over-arc points to the under-arc that receives the product. Moreover this quandle is invariant under the Reidemeister moves thus constituting a knot invariant called the knot quandle. Upon fixing a finite quandle $X$ the number of homomorphisms (also known as colorings) from the knot quandle to the finite quandle $X$ is also a knot invariant. This invariant is referred to as counting colorings. An assignement of elements of $X$ (colors) to the knot quandle is a homomorphism if at the crossings the colors assigned to the arcs satisfy the same relations as the arcs do. One efective way of finding colorings of a knot by a given finite quandle $X$ resorts to braids (\cite{Birman}). Any knot can be represented by the closure of a certain braid so we assign colors to the strands at the top of the braid (the input colors) whose closure yields the knot under study. These colors propagate down the braid giving rise to knew colors as one goes past each crossing and eventually obtaining a new string of colors on the strands at the bottom of the braid (the output colors). If this assignment is a coloring then the string of colors at the top equals the string of colors at the bottom, otherwise it is not a coloring. In some way, the braid acts on the input colors to produce the output colors; if the output colors match the input colors then these input colors stand for a coloring, otherwise they do not stand for a coloring. We will assume that the finite labeling quandles we will be using are all Alexander quandles. Furthermore, we use the Burau representation (\cite{Birman}) to perform the calculations. The braid group on $N$ strands, $B\sb{N}$, is generated by the so-called standard generators, $\sigma\sb{i}$, $i=1, \dots , N-1$ (\cite{Birman}). Geometrically, $\sigma\sb{i}$ is the $(i+1)$-th strand going over the $i$-th strand. In the Burau representation the standard generators are mapped to matrices over the ring of Laurent polynomials in the variable $T$ as follows:
\[
\sigma_i \longmapsto I_{i-1}\otimes \left[ \begin{array}{cc} 0 & T \\ 1 &
1-T \end{array}\right] \otimes I_{N-1-i}
\]
Then the action of the braid on the input color plus the matching against the input we referred to, translates into the ``eigenvector'' equation
\[
( \, a\sb{1} \, a\sb{2} \, \dots \, a\sb{N}\, ) B (b) = ( \, a\sb{1} \, a\sb{2} \, \dots \, a\sb{N}\, )
\]
where $a\sb{i}, i=1, \dots , N$ are the colors from the Alexander quandle $X$ and $B(b)$ is the Burau matrix of the braid $b=\sigma\sb{i\sb{1}}\sp{\epsilon\sb{i\sb{1}}}\cdot \dots \cdot \sigma\sb{i\sb{L}}\sp{\epsilon\sb{i\sb{L}}}$ whose closure yields the knot under study. If the equality holds then $(a\sb{1}, \dots , a\sb{N})$ stands for a coloring, otherwise it does not stand for a coloring. Moreover, note that the equality is understood in the quotient corresponding to the Alexander quandle we are dealing with.

\bigbreak

As sets, these Alexander quandles are finite. Moreover, they have a ring structure induced by addition and multiplication of polynomials. For a given positive integer $N$, the ring of $N\times N$ matrices with entries from one of these Alexander quandles is thus also finite. Since each Burau matrix, say $B$, is invertible, then it has finite order i.e., there is a (positive) integer $M$ such that $B\sp{M}$ is the identity matrix.

\bigbreak

Given a knot $K\sb{1}$ and a braid $b\in B\sb{N}$ whose closure yields $K\sb{1}$, consider the sequence of knots given by $K\sb{n}=\widehat{b\sp{n}}$, where $\widehat{\dots }$ denotes ``closure''. When coloring the knots from this sequence with an Alexander quandle, we write down the ``eigenvector'' equation as above:
\[
( \, a\sb{1} \, a\sb{2} \, \dots \, a\sb{N}\, ) \bigl( B (b)\bigr)\sp{n} = ( \, a\sb{1} \, a\sb{2} \, \dots \, a\sb{N}\, )
\]
for knot $K\sb{n}$ and where $B(b)$ is the Burau matrix for $b$ in the appropriate quotient. Since there is an $M$ such that $B(b)\sp{M}$ is the identity matrix, then there is periodicity in the colorings of the knots of this sequence, in particular, in the sense that there are only $M$ distinct sets of eigenvectors for this sequence of knots. Moreover there are at most $M$ distinct distributions of colors per arcs over any $M$ consecutive factors of $b\sp{n}$. Also, for a fixed $b\sp{n}$, suppose $n=lM+r$; there are $l$ identical coloring patterns followed by the coloring pattern associated to the $b\sp{r}$.

This periodicity issue was first discussed in \cite{jsCetal0} if only for some sequences of torus knots. We are now claiming it in more generality, both for all Alexander quandles and not only for the particular case of sequences of torus knots.

\bigbreak

Now assume that an abelian group $A$ and a $2$-co-cycle $\phi$ have been chosen besides the Alexander quandle $X$. Let $|A|$ stand for the cardinality of $A$ and let $M$ be as above i.e., $B(b)\sp{M}$ is the identity matrix. Given a positive integer $n$ let
\[
n=lM|A|+r \qquad \qquad \text{ with } 0\leq r < M|A|
\]
We let $c(b\sp{m})$ denote the set of crossings of $b\sp{m}$ and $1\sb{A}$ denote the identity in the group $A$. The CJKLS invariant of $K\sb{n}=\widehat{b\sp{n}}$ then reads
\begin{align*}
Z(\widehat{b\sp{n}})&=\sum\sb{\underset{\text{s.t. ...}}{a\sb{1}, \dots , a\sb{N} \in X}}  \prod\sb{\tau\in c(b\sp{n}) }\phi\sp{\epsilon\sb{\tau}} = \sum\sb{\underset{\text{s.t. ...}}{a\sb{1}, \dots , a\sb{N} \in X}}\Biggl(  \prod\sb{\tau\in c(b\sp{M})}\phi\sp{\epsilon\sb{\tau}}  \Biggr) \sp{l|A|}\cdot  \prod\sb{\tau\in c(b\sp{r})}\phi\sp{\epsilon\sb{\tau}}=\\
&=\sum\sb{\underset{\text{s.t. ...}}{a\sb{1}, \dots , a\sb{N} \in X}}\Biggl(  \Biggl( \prod\sb{\tau\in c(b\sp{M})}\phi\sp{\epsilon\sb{\tau}}  \Biggr) \sp{|A|}\Biggr)\sp{l}\cdot  \prod\sb{\tau\in c(b\sp{r})}\phi\sp{\epsilon\sb{\tau}}=\sum\sb{\underset{\text{s.t. ...}}{a\sb{1}, \dots , a\sb{N} \in X}}\bigl(  1\sb{A}\bigr)\sp{l}\cdot  \prod\sb{\tau\in c(b\sp{r})}\phi\sp{\epsilon\sb{\tau}}=\\
&=\sum\sb{\underset{\text{s.t. ...}}{a\sb{1}, \dots , a\sb{N} \in X}}\prod\sb{\tau\in c(b\sp{r})}\phi\sp{\epsilon\sb{\tau}}
\end{align*}
The dots next to s.t. and under $a\sb{1}, \dots , a\sb{N} \in X$ stand for the conditions the $a\sb{i}$ have to satisfy in order to stand for colorings; these conditions stem from the eigenvector equation and there are at most $M$ distinct conditions. Since $0 \leq r < M|A|$ then there are at most $M|A|$ different values of the CJKLS invariant for this sequence of knots $K\sb{n}$. Since we regard the CJKLS invariant as taking on values in the euclidean metric space $\mathbb{R}\sp{|A|}$ then the sequence of the CJKLS invariant for $K\sb{n}$ is bounded and so is the sequence of the free energies for this sequence of knots. If the crossing number of $K\sb{n}$ goes to infinity as $n$ goes to infinity then the free energy per crossing of this sequence of knots goes to the null vector.

\bigbreak

\section{Example}

Fix a positive integer $N$ and consider the sequence of torus knots $\bigl( T(N, n)\bigl) \sb{n\in \mathbb{N}}$. We have
\[
T(N, n) = \biggl( \sigma\sb{N-1}\sigma\sb{N-2}\cdots \sigma\sb{2}\sigma\sb{1}\biggr)\sp{n}
\]
Moreover the crossing number of $T(N, n)$ is
\[
\min \{ |N|(|n|-1), |n|(|N|-1)  \}
\]
(\cite{Murasugi}) which tends to infinity as $n$ tends to infinity. Thus, the sequence of the free energies per crossing from a CJKLS invariant with an Alexander quandle as labeling quandle, for a sequence of torus knots as the one indicated tends to the null vector thereby giving rise to a hyperfinite knot.

\begin{theacknowledgments}
The author acknowledges finantial support from POCTI/FCT/FEDER.
\end{theacknowledgments}


\end{document}